\documentclass[graybox]{svmult}

\usepackage{mathptmx}       
\usepackage{helvet}         
\usepackage{courier}        
\usepackage{type1cm}        
\usepackage{makeidx}         
\usepackage{graphicx}        
                  
\usepackage{multicol}        
\usepackage[bottom]{footmisc}
 \usepackage{float}
\usepackage{url}
\usepackage{cite}
\usepackage{mathrsfs} 

\usepackage{lineno}

\usepackage{color}
\usepackage{array}
\usepackage{bigstrut}
\usepackage{epsfig}
\usepackage{amssymb}
\usepackage{amsmath}
\usepackage{amstext,amsfonts} 		
\usepackage{latexsym,amscd,enumerate,supertabular}

\usepackage{hyperref}
\hypersetup{
    colorlinks,
    hidelinks,
    citecolor=blue,
    filecolor=black,
    linkcolor=red,
    urlcolor=green,
    plainpages=false,
    pdfpagelabels=true
}

\newcommand{\q}{\quad}
\newcommand{\ee}{{\rm e}\hspace{1pt}}
\spnewtheorem{assumption}{Assumption}{\bfseries}{\rmfamily}
\newcommand{\dd}{\hspace{0.5pt}{\rm d}\hspace{0.5pt}}
\newcommand{\mat}[1] {{\mathbf{#1}}}

\makeindex             

\begin{document}

\title*{Explicit Exponential Rosenbrock Methods and their Application in Visual Computing%
\thanks{This work has been partially supported by King Abdullah University of Science and Technology (KAUST baseline funding).}
}
\author{Vu Thai Luan and Dominik L. Michels}
\institute{Vu Thai Luan \at  Department of Mathematics, Southern Methodist University,
  PO Box 750156, Dallas, TX 75275-0156, USA, \email{vluan@smu.edu} 
\and {Dominik L. Michels \at  ​​​​Computational Sciences Group, Visual Computing Center, King Abdullah University of Science and Technology, Thuwal, 23955, KSA, \email{dominik.michels@kaust.edu.sa}}}

\maketitle
\abstract{We introduce a class of explicit exponential Rosenbrock methods for the time integration of large systems of stiff differential equations. Their application with respect to simulation tasks in the field of visual computing is discussed where these time integrators have shown to be very competitive compared to standard techniques. In particular, we address the simulation of elastic and nonelastic deformations as well as collision scenarios focusing on relevant aspects like stability and energy conservation, large stiffnesses, high fidelity and visual accuracy.}


\section{Introduction}
\label{sec:introduction}
Developing numerical models for practical simulations in science and engineering usually results in problems regarding the presence of wide-range time scales. These problems involve both slow and fast components leading to rapid variations in the solution. This gives rise to the so-called \emph{stiffness phenomena}. Typical examples are models in molecular dynamics (see e.g.~\cite{Michels:2015}), chemical kinetics, combustion, mechanical vibrations (mass-spring-damper models), visual computing (specially in computer animation),  computational fluid dynamics,  meteorology, etc., just to name a few. They are usually formulated as systems of stiff differential equations which can be cast in the general form 
\begin{equation} \label{eq1}
u'(t)=F(u(t)), \q  u(t_0)=u_0,
\end{equation}
where $u\in \mathbb{R}^n$ is the state vector and $F: \mathbb{R}^n \longrightarrow \mathbb{R}^n$ represents the vector field. The challenges in solving this system are due to its stiffness by means of the eigenvalues of the Jacobian matrix of $F$  differing by several orders of magnitude. In the early days of developing numerical methods for ordinary differential equations (ODEs), classical methods such as the explicit Runge--Kutta integrators were proposed.  For stiff problems, however, they are usually limited by stability issues due to the CFL condition leading to the use of unreasonable time steps, particularly for large-scale applications. The introduction of implicit methods such as semi-implicit, IMEX (see \cite{Ascher97}), and BDF methods (see \cite{Curtiss1952,Gear1971}) has changed the situation. Theses standard methods require the solution of nonlinear systems of equations in each step. As the stiffness of the problem increases, considerably computational effort is observed. This can be seen as a shortcoming of the implicit schemes.
 
 In the last twenty years, with the new developments of numerical linear algebra algorithms in computing matrix functions \cite{Hochbruck98,AH11,NW12}, exponential integrators have become an alternative approach for stiff problems (see the survey 
\cite{HO10}; next to physics simulations, exponential integrators are nowadays also employed for different applications as for the construction of hybrid Monte Carlo algorithms, see \cite{pmlr-v37-chao15}). For the fully nonlinear stiff system \eqref{eq1}, we mention good candidates, the so-called \emph{explicit exponential Rosenbrock methods}, which can handle the stiffness of the system in an explicit and very accurate way. This class of exponential integrators was originally proposed in \cite{HO06} and further developed in \cite{HOS09,LO14a,LO16,LO13}. They have shown to be very efficient both in terms of accuracy and computational savings. In particular, the lower-order schemes were recently successfully applied to a number of different applications \cite{Gondal2010,Geiger2012,Tambue2013,Zhuang2015,Chen2017} and very recently the fourth- and fifth-order schemes were shown to be the method of choice for some meteorological models (see \cite{Luan18}).

In this work, we show how the exponential Rosenbrock methods (particularly higher-order schemes) can be also applied efficiently in order to solve problems in computational modeling of elastodynamic systems of coupled oscillators (particle systems) which are often used in visual computing (e.g.~for computer animation). In their simplest formulation, their dynamics can be described using Newton's second law of motion leading to a system of second-order ODEs of the form
\begin{equation} \label{eq2}
m_i \ddot{x}_i + \sum_{ j \in \mathcal{N}(i)} k_{ij} (\| x_i-x_j\|- \ell_{ij}) \frac{x_i-x_j}{\| x_i-x_j\|} = g_i (x_i, \dot{x}_i,\cdot), \q i=1,2,\cdots,N,
\end{equation} 
where $N$ is the number of particles, $x_i \in \mathbb{R}^3$, $m_i$, $k_{ij} $, $\ell_{ij}$ denote the position of particle $i$ from the initial position, its mass,  the spring stiffness, the equilibrium length of the spring between particles $i$ and  $j$, respectively, and $\mathcal{N}(i)$ denotes the set of indices of particles that are connected to particle $i$ with a spring (the neighborhood of particle  $i$). Finally, $g_i$ represents the external force acting on particle $i$ which can result from an external potential, collisions, etc., and can be dependent of all particle positions, velocities, or external forces set by user interaction.

Our approach for integrating \eqref{eq2}  is first to reformulate it in the form of \eqref{eq1} (following a novel approach in \cite{Michels2017}).   The reformulated system is a very stiff one since the linear spring forces usually  possess very high frequencies. Due to the special structure of its linear part (skew-symmetric matrix) and large nonlinearities, we then make use of exponential Rosenbrock methods. Moreover, we propose to use the improved algorithm in \cite{Luan18} for the evaluation of a linear combination of $\varphi$-functions acting on certain vectors $v_0,\ldots,v_p$, i.e. $\sum_{k=0}^{p}  \varphi_k (A)v_k $
which is crucial for implementing exponential schemes. 
Our numerical results on a number of complex models in visual computing
indicate that this approach significantly reduces computational time over the current state-of-the-art techniques while maintaining sufficient levels of accuracy.

This contribution is organized as follows. In Section~\ref{sec:coupled_oscillators}, we present a reformulation of systems of coupled oscillators \eqref{eq2} in the form of \eqref{eq1} and briefly review previous approaches used for simulating these systems in visual computing.  In Section~\ref{sec:exponential_methods}, we describe the exponential Rosenbrock methods as an alternative approach for solving large stiff systems \eqref{eq1}. The implementation of these methods is discussed in Section~\ref{sec:implement}, where we also introduce a new procedure to further improve one of the state-of-the-art algorithms.
In Section~\ref{sec:experiments} we demonstrate the efficiency of the exponential Rosenbrock methods on a number of  complex models in visual computing. In particular, we address the simulation of deformable bodies, fibers including elastic collisions, and crash scenarios including nonelastic deformations. These examples focus on relevant aspects in the realm of visual computing, like stability and energy conservation, large stiffness values, and high fidelity and visual accuracy. We include an evaluation against classical and state-of-the-art methods used in this field.
Finally, some concluding remarks are given in Section~\ref{sec:conclusion}.

\section{Reformulation of Systems of Coupled Oscillators}
\label{sec:coupled_oscillators}
We first consider the system of coupled oscillators \eqref{eq2}. Let
$x(t)\in\mathbb{R}^{3N}$, $M\in \mathbb{R}^{3N\times 3N} $,  $D\in \mathbb{R}^{3N\times 3N} $, $K\in \mathbb{R}^{3N\times 3N} $ and $g(x) \in \mathbb{R}^{3N}$  denote the vector of positions, 
the mass matrix (often diagonal and thus nonsingular), the damping matrix, the spring matrix (stiff), and the total external forces acting on the system, respectively. Using these matrix notations and  denoting  $A=M^{-1}K$, \eqref{eq2} can be written as a system of second-order ODEs
\begin{equation} \label{eq2.1}
x''(t)+A x(t)= g(x(t)), \q  x(t_0)=x_0, \ x'(t_0)=v_0. 
\end{equation}
Here $x_0, v_0$ are some given initial positions and velocities.  For simplicity we neglect damping and 
 assume that $A$ is a symmetric, positive definite matrix (this is a reasonable assumption in many models, see \cite{Michels:2014}). Therefore, there exists a unique positive definite matrix $\Omega$ such that $A=\Omega^2$ (and clearly $\Omega^{-1}$ exists).

Following our approach in \cite{Michels2017}, we introduce the new variable 
\begin{equation} \label{eq2.2}
u(t)=\left[ 
\begin{array}{c}
    \Omega x(t)    \\
    x'(t)   
\end{array} \right].
\end{equation}
Using this one can reformulate \eqref{eq2.1} as a first-order system of ODEs of the form like \eqref{eq1}:
\begin{equation} \label{eq2.3}
u'(t)=F(u(t))= \mathscr{A}u(t)+G(u(t)), \q u(t_0)=u_0,
\end{equation}
where
\begin{equation} \label{eq2.4}
\quad\mathscr{A}=
\left[ \begin{array}{ccc}
    {\bf 0} & \Omega  \\
    -\Omega & {\bf 0}
\end{array} \right], \q
G(u)=
\left[ \begin{array}{c}
    {\bf 0}   \\
    g(x)  
\end{array} \right].
\end{equation}
Since the linear spring forces usually possess high frequencies (thus $\|K\|\gg 1$ and so is $\|A\|$),  \eqref{eq2.3} becomes a very stiff ODE. 
Regarding the new formulation \eqref{eq2.3}-\eqref{eq2.4}, we observe the following two remarks.

\begin{remark}\label{remark1}
 Clearly, the  linear part associated with $\mathscr{A}$ is a skew-symmetric matrix. We note that this significantly differs from the common way of reformulating  \eqref{eq2.1} that is to use the change of variable 
$
 X(t)=[ x(t) , \  x'(t)]^{T}   
$ which results in a non-symmetric matrix. The great advantage of \eqref{eq2.4} is that we know the nonzero eigenvalues of $\mathscr{A}$ are all pure imaginary and are in pairs $\pm \lambda_k i$. Moreover, one realizes that  $\mathscr{A}$ is an infinitesimal symplectic (or Hamiltonian) since $$J\mathscr{A}_{\text{new}}=\left[ \begin{array}{ccc}
   -\Omega & {\bf 0}\\
   {\bf 0} &  -\Omega
\end{array} \right]$$ is symmetric, i.e., $ J\mathscr{A}=(J \mathscr{A})^T$, where $$J=\left[ \begin{array}{ccc}
    {\bf 0} & I  \\
    -I & {\bf 0}
\end{array} \right]\,.$$
\end{remark}
\begin{remark}\label{remark2}
If the Jacobian matrix $F'(u)=\mathscr{A}+ G'(u)$ is infinitesimal symplectic, \eqref{eq2.3} is a Hamiltonian system. This can be fulfilled since a typical situation in Hamiltonian systems is that $g(x)=\nabla f(x)$ for some function $f(x)$ and thus $g'(x)= \nabla^2 f(x)$ becomes a Hessian matrix, which is symmetric. 
\end{remark}

As seen, either using the common way (mentioned in Remark~\ref{remark1}) or the new way  \eqref{eq2.2} for reformulating  \eqref{eq2.1}, one has to solve the stiff ODE \eqref{eq2.3}. In visual computing it is usually solved by explicit methods such as the fourth-order Runge--Kutta methods, semi-implicit methods such as  the St\"ormer--Verlet methods, the backward differentiation formulas (BDF-1 and BDF-2) methods, or IMEX methods. In this regard, we refer to some contributions in the context of interacting deformable bodies,  cloth, solids, and elastic rods, see  \cite{Terzopoulos:1987a,Baraff:1998,Eberhardt:2000,Hauth:2001,Goldenthal:2007,Bergou:2008}. For large-scale applications associated with stiff systems, however, both types of these time integration techniques have their own limitations as mentioned in the introduction. In recent years, exponential integrators have
shown to be competitive for large-scale problems in physics and for nonlinear parabolic PDEs, as well as for highly oscillatory problems (see \cite{HO10}). They have attracted much attention by the broad computational mathematics community since mid-1990s. At the time while solving linear systems $(I-\alpha h J)x=v$ with some Jacobian matrix $J$  (required when using implicit methods)  is generally only linear convergence, it was realized that Krylov subspace methods for approximating the action of a matrix exponential on a vector, $\ee^{hJ} v$, offer superlinear convergence (see \cite{HL97}). Unless a good preconditioner is available, this is clearly a computational advantage of exponential integrators over implicit methods. This has been addressed in the visual computing community very recently through a number of interesting work on exponential integrators, see e.g.\cite{Michels:2014,Michels:2015b,Michels:2016:DCM:2988458.2988464,Michels2017}. Inspired by this interest, in the following sections we will show how exponential Rosenbrock methods -- one of the popular classes of exponential integrators -- can be applied for simulating systems of coupled oscillators. 

\section{Explicit Exponential Rosenbrock Methods}
\label{sec:exponential_methods}
In this section, based on \cite{HO06,HOS09,LO14a,LO16,L17} we present a compact summary of the introduction of exponential Rosenbrock methods and their derivations for methods of order up to 5. We then display some efficient schemes for our numerical experiments for some applications in visual computing. 
\subsection{Approach}
Motivated by the idea of deriving Rosenbrock-type methods, see \cite[Chap. IV.7]{HW96}, instead of integrating the fully nonlinear system \eqref{eq1} (which has a large nonlinearity for stiff problems), one can replace it by a sequence of semilinear problems 
\begin{equation} \label{eq3.1}
u'(t)=F(u(t))=J_n u(t)+g_n(u(t)),
\end{equation}
by linearizing the forcing term $F(u)$ in each time step at the numerical solution $u_n$ (due to \cite{Pope63}) with
\begin{equation} \label{eq3.2}
J_n =F'(u_n), \ g_n(u)=F(u)-J_n u
\end{equation}
are the Jacobian and the nonlinear remainder, respectively. An advantage of this approach is that $g'_n(u_n)=F'(u_n)-J_n=0$ which shows that the new nonlinearity $g_n (u)$ has a much smaller Lipschitz constant than that of the original one $F(u)$.
The next idea is to handle the stiffness by solving the linear part $J_n u$ exactly and integrating the new nonlinearity $g_n (u)$  explicitly.  
For that, the representation of the exact solution at time $t_{n+1}=t_n+h$ of \eqref{eq3.1} using the variation-of-constants formula 
\begin{equation} \label{eq3.3}
u(t_{n+1})=\ee^{hJ_n}u(t_n) +\int_{0}^{h} \ee^{(h-\tau)J_n} g_n(u(t_n+\tau)) \dd \tau
\end{equation}
plays a crucial role in constructing this type of integrators. As seen from \eqref{eq3.3}, while the linear part can be integrated exactly by computing the action of the matrix exponential $\ee^{hJ_n}$ on the vector $u(t_n)$, the integral involving  $g_n (u)$ can be approximated by some quadrature. This procedure results in the so-called \emph{exponential Rosenbrock methods}, see \cite{HO06,HOS09}.

\begin{remark}
\label{remark3.1}
\textnormal{
For the system of coupled oscillators \eqref{eq2}, the forcing term $F(u)$ has the semilinear form  \eqref{eq2.3}, which can be considered as a fixed linearization problem \eqref{eq3.1} (i.e. $J_n=\mathscr{A}$). Therefore, one can directly apply explicit the exponential Runge--Kutta methods (see \cite{HO05}) to  \eqref{eq2.3}. The advantage of these methods is that the time-step $h$ is not restricted by the CFL condition when integrating the linear part $\mathscr{A}u$. In our applications, however, the nonlinearity $G(u)$ is  large in which the CFL condition usually serves as a reference for setting the time-step. In particular, for the stability $hL_{G}$ should be sufficiently small  ($L_{G}$ is the Lipschitz constant of $G(u)$). In this regard,  the dynamic linearization approach  \eqref{eq3.1} applied to  \eqref{eq2.3}  
\begin{equation} \label{eqRemark3.1}
u'(t)=F(u)= \mathscr{A}u +G(u)=J_n u+G_n(u)
\end{equation}
 with 
\begin{equation} \label{eqJn}
J_n = \mathscr{A}+G'(u_n),
\end{equation}
 offers a great advantage in improving the stability (in each step) when integrating $G(u)$. This is because instead of integrating the original semilinear problem with  large nonlinearity $G(u)$,
we only have to deal with a much smaller nonlinearity $G_n(u)$   (as mentioned above). Note that the new linear part $J_n u$ with the Jacobian $J_n$ now incorporates both $\mathscr{A}$ and the Jacobian of the nonlinearity $G(u)$, which can be again solved exactly. It is thus anticipated that this idea of exponential Rosenbrock methods opens up the possibility to take even larger time steps compared to exponential Runge--Kutta methods.
}
\end{remark}

\subsection{Formulation of a Second-order and General Schemes}
In this subsection, we will illustrate the approach of exponential Rosenbrock methods by presenting a simple derivation of a second-order scheme and formulating general schemes.  
\subsubsection{A Second-order Scheme}
First, expanding  $u(t_n+\tau)$ in a Taylor series gives
$u(t_n+\tau)=u(t_n)+\tau u'(t_n)+ \mathcal{O}(\tau^2)$. Then inserting this into  $g_n(u(t_n+\tau))$ and again expanding it as a Taylor series around $u(t_n)$ (using $g'_n(u(t_n))=0$) leads to
\begin{equation} \label{eq3.4}
g_n (u(t_n+\tau))=g_n (u(t_n)) +\mathcal{O}(\tau^2)\,.
\end{equation}
Inserting \eqref{eq3.4} into  the integral part of \eqref{eq3.3} and denoting $\varphi_1 (h J_n)=\frac{1}{h} \int_{0}^{h} \ee^{(h-\tau) J_n} \dd \tau$ gives
\begin{equation}\label{eq3.4a}
u(t_{n+1})=\ee^{hJ_n} u(t_n) +h \varphi_1 (h J_n) g_n(u(t_n))  +\mathcal{O}(h^3).
\end{equation}
 Neglecting the local error term $\mathcal{O}(h^3)$ results in a second-order scheme, which can be reformulated as 
\begin{equation} \label{eq3.5}
u_{n+1}=u_n + h\varphi_1(hJ_n)F(u_n)
\end{equation}
by replacing $g_n(u(t_n))$ by  \eqref{eq3.2} and using the fact that $\varphi_1 (z)=(e^z-1)/z$. This scheme was derived before and named as \emph{exponential Rosenbrock-Euler method}, see \cite{HO06,HOS09} (since
when considering the formal limit  $J_n\rightarrow \mathbf{0}$, \eqref{eq3.5} is the underlying Euler method). The derivation here, however, shows directly that this scheme has an order of consistency three and thus it is a second-order stiffly accurate method (since the constant behind the Landau notation $\mathcal{O}$ only depends on the regularity assumptions on $u(t)$ and $g_n(u)$, but is independent of $\|J_n\|$).
\subsubsection{General Schemes}
For the derivation of higher-order schemes, one can proceed in a similar way as the construction of classical Runge--Kutta methods. 
Namely, one can approximate the integral in \eqref{eq3.3} by using some higher-order quadrature rule with nodes $c_i$ in $[0,1]$ and weights $b_i(h J_n)$ which are matrix functions of $hJ_n$, yielding  
\begin{equation} \label{eq3.8}
u(t_{n+1})\approx \ee^{h J_n}u(t_n) + h \sum_{i=1}^{s} b_i(h J_n)g_n( u(t_n+c_i h)).
\end{equation}
The unknown intermediate values  $u(t_n+c_i h)$ can be again approximated by using  \eqref{eq3.3} (with $c_i h$ in place of $h$) with another quadrature rule using the same nodes $c_j$, $1\leq j \leq i-1$, (to avoid generating new unknowns) and new weights $a_{ij}(hJ_n)$, leading to
\begin{equation} \label{eq3.9}
u(t_n+c_i h)\approx \ee^{c_i h J_n}u(t_n) + h_n\sum_{j=1}^{i-1} a_{ij}(h J_n)g_n(u(t_n+c_j h)).
\end{equation}
Let us denote  $u_n \approx u(t_n)$ and $U_{ni}\approx u(t_n +c_i h_n)$.  As done for \eqref{eq3.5}, using  \eqref{eq3.4} (with $c_i h, h$ in place of $\tau$, respectively) one can reformulate \eqref{eq3.8} and \eqref{eq3.9} in a similar manner, which yields the general format of $s$-stage explicit exponential Rosenbrock methods 
\begin{subequations} \label{eq3.10}
\begin{align}
U_{ni}&= u_n + c_i h \varphi _{1} ( c_i h J_n)F(u_n) +h \sum_{j=2}^{i-1}a_{ij}(h J_n) D_{nj}, \label{eq3.10a} \\
u_{n+1}& = u_n + h\varphi _{1} ( h J_n)F(u_n) + h \sum_{i=2}^{s}b_{i}(h J_n) D_{ni}  \label{eq3.10b}
\end{align}
with
\begin{equation} \label{eqDni}
  D_{ni}= g_n ( U_{ni})- g_n(u_n ), 
\end{equation}
\end{subequations}
As in  \eqref{eq3.4}, we have $D_{ni}=\mathcal{O}(h^2)$. Thus, the general methods \eqref{eq3.10} are small perturbations of the exponential Rosenbrock-Euler method  \eqref{eq3.5}.
Note that the weights $a_{ij}(hJ_n)$ and $b_i(h J_n)$ are usually linear combinations of  $\varphi _{k} (c_i h J_n)$ and $\varphi_{k} (h J_n),$ respectively, where the
 $\varphi$ functions (similar to $\varphi_1$) are given by
\begin{equation} \label{eq3.6}
\varphi _{k}(hZ)=\frac{1}{h^k}\int_{0}^{h} \ee^{(h-\tau )Z} \tau^{k-1} \dd \tau , \quad k\geq 1
\end{equation}
and satisfy the recursion relation
\begin{equation} \label{eq3.7}
\varphi _{k+1}(z)=\frac{\varphi _{k}(z)-\frac{1}{k!}}{z}, \quad k\geq 1.
\end{equation}
It is important to note that these functions are bounded (uniformly) independently of $\|J_n\|$ (i.e. the stiffness) so do the coefficients  $a_{ij}(hJ_n)$ and $b_i(h J_n)$ (see e.g. \cite{HO10}). 

Clearly, using exponential Rosenbrock schemes  \eqref{eq3.10} offers some good advantages. First, they are fully explicit and do not require the solution of linear or nonlinear systems of equations. Second, as mentioned above, they offer a better stability when solving stiff problems with large nonlinearities and thus allow to use larger time-steps. Third, since the Jacobian of the new nonlinearity vanishes at every step ($g'_n(u_n)=0$), the derivation of the order conditions and hence the schemes can be simplified considerably. In particular, higher-order stiffly accurate schemes can be derived with only a few stages (see the next section).  

The convergence analysis of exponential Rosenbrock methods is usually carried out in an appropriate framework (strongly continuous semigroup) under regularity assumptions on the solution $u(t)$ (sufficiently smooth) and $g_n(u)$ (sufficiently Fr\'echet differentiable in a neighborhood of the solution) with uniformly bounded derivatives in some Banach space. For more details, we refer to \cite{HOS09,LO14a}.

\subsection{Selected Schemes for Numerical Simulations}
\label{sec:selected_exponential_schemes}
First, we discuss some important points for the derivation of  exponential Rosenbrock schemes. Clearly, the unknown coefficients $a_{ij}(hJ_n)$ and $b_i(h J_n)$ has to be determined by solving order conditions. For nonstiff problems, where the Jacobian matrix has a small norm, one can expand those matrix functions using classical Taylor series expansions, leading to nonstiff order conditions and in turn classical exponential schemes (see e.g. \cite{Cox2002,Krogstad2005}). For stiff problems, however, one has to be cautious when analyzing the local error to make sure that error terms do not involve powers of $J_n$ (which has a large norm). Recently, Luan and Ostermann \cite{LO12b,LO13}  derived a new expansion of the local error which fulfills this requirement and thus derived  a new stiff order conditions theory for methods of arbitrary order (both for exponential Runge--Kutta and exponential Rosenbrock methods). As expected, with the same order, the number of order conditions for exponential Rosenbrock methods is significant less than those for  exponential Runge--Kutta methods. For example, in  Table~\ref{tb3.1}, we display the required 4 conditions for deriving schemes up to order 5 in \cite{LO14a} (note that  for exponential Runge--Kutta methods, 16 order conditions are required for deriving schemes of order 5, see \cite{LO14b}). 
\renewcommand{\arraystretch}{1.4}%
\begin{table}[H]
\caption{Stiff order conditions for exponential Rosenbrock methods up to order five. Here $Z$ and $K$ denote arbitrary square matrices and $\psi_{3,i}(z)= \sum_{k=2}^{i-1}a_{ik}(z)\frac{c^2_k}{2!}-c^{3}_i \varphi_{3} ( c_i z)$. }
\begin{center}
\begin{tabular}{ |c|c|c| }
\hline
\textbf{No.} & \textbf{Order Condition} & \textbf{Order} \\
\hline
1&$\sum_{i=2}^{s} b_i (Z)c^2_i=2\varphi_3 (Z) $&3 \\
 \hline
2&$\sum_{i=2}^{s} b_i (Z)c^3_i=6\varphi_4 (Z) $&4 \\
\hline
3&$\sum_{i=2}^{s} b_i (Z)c^4_i=24\varphi_5 (Z) $&5 \\
4&$\sum_{i=2}^{s} b_i (Z)c_i  K \psi _{3,i}(Z)=0 $&5 \\
\hline
\end{tabular}
\end{center}
\label{tb3.1}
\end{table}%
\vspace{-10pt}
We note that with these order conditions one can easily derive numerous different schemes of order up to 5. Taking the compromise between efficiency and accuracy into consideration, we seek for the most efficient schemes for our applications. Namely, the following two representative fourth-order schemes are selected. \\
\\
$\mathtt{exprb42}$ (a fourth-order 2-stage scheme which can be  considered as a superconvergent scheme, see \cite{L17}):
\begin{subequations} \label{eq3.11}
\begin{align}
 U_{n2}&= u_n + \tfrac{3}{4} h \varphi _{1} ( \tfrac{3}{4} h J_n)F(u_n), \label{eq3.11a} \\
u_{n+1}& = u_n + h \varphi _{1} ( h J_n)F(u_n) + h \tfrac{32}{9}\varphi_3 (h J_n) (g_n (U_{n2})- g_n(u_n )) \label{eq3.11b}.  
\end{align}
\end{subequations}
\\
$\mathtt{pexprb43}$ (a  fourth-order 3-stage scheme, which can be implemented in parallel, see  \cite{LO16}):
\begin{subequations} \label{eq3.12}
\begin{align}
U_{n2}&= u_n + \tfrac{1}{2} h \varphi _{1} ( \tfrac{1}{2} h J_n)F(u_n), \\
U_{n3}&= u_n +  h \varphi _{1} (h J_n)F(u_n), \\
u_{n+1}& = u_n + h \varphi _{1} ( h J_n)F(u_n) + h \varphi_3 (h J_n) (16 D_{n2}-2 D_{n3})\nonumber\\
& \,\,\,\,\,\,\,\,\,\,\,\,\,+ h \varphi_4 (h J_n) (-48 D_{n2}+12 D_{n3}).
\end{align}
\end{subequations}
Note that  the vectors $D_{n2}$ and $D_{n3}$ in \eqref{eq3.12} are given by \eqref{eqDni}, i.e., 
 $D_{n2}=g_n (U_{n2})- g_n(u_n )$ and $D_{n3}=g_n (U_{n3})- g_n(u_n )$.

\section{Implementation}
\label{sec:implement}
In this section, we present the implementation of exponential Rosenbrock methods for the new formulation \eqref{eq2.3} of the system of coupled oscillators. First, we discuss on the computation of the matrix square root $\Omega$ needed for the reformulation. We then briefly review some state-of-the-art algorithms for implementing exponential Rosenbrock methods and introduce a new routine which is an improved version of one of these algorithms (proposed very recently in \cite{Luan18}) for achieving more efficiently. Finally, we specifically discuss applying this routine for implementing the selected schemes $\mathtt{exprb42}$ and $\mathtt{pexprb43}$. 

\subsection{Computation of the Matrix Square Root $\Omega=\sqrt{A}$}
For the computation of $\Omega=\sqrt{A}$ used in  \eqref{eq2.3}, we follow our approach in \cite{Michels2017}. Specifically, we use the Schur decomposition for moderate systems. For large systems, the Newton square root iteration (see \cite{Higham08}) is employed in order to avoid an explicit precomputation of $\Omega$. Namely, one can use the following simplified iteration method for approximating the solution of  the equation $\Omega^2= A$:
 \begin{itemize}
\item [(i)] choose $\Omega_0=A$ ($k=0$),
\item [(ii)]  \ update  $\Omega_{k+1}= \frac{1}{2}(\Omega_k+\Omega^{-1}_k A)$.
\end{itemize}
This method offers unconditional quadratic convergence with much less cost compared to the Schur decomposition. We note that  $\Omega^{-1}$ can be computed efficiently using a Cholesky decomposition since $\Omega$ is symmetric and positive definite and it is given by $\Omega^{-1}=\mat{S}^{-1}\mat{S}^{-\mathsf{T}}$, where  $\mat{S}$ is an upper triangular matrix with real and positive diagonal entries. For more details, we refer to  \cite{Higham08,Michels2017}.

With $\Omega$ at hand, one can easily compute the Jacobian $J_n$ as in \eqref{eqJn} and $F(u), G_n(u)$ as in \eqref{eqRemark3.1}. As the next step, we  discuss the implementation of the exponential Rosenbrock schemes. 
\subsection{Implementation of Exponential Rosenbrock Methods}
In view of the exponential Rosenbrock schemes in Section~\ref{sec:exponential_methods}, each stage requires  the evaluation of a linear combination of $\varphi$-functions acting on certain vectors $v_0,\ldots,v_p$
\begin{equation} \label{eq3.14}
\varphi_0 (M)v_0 + \varphi_1 (M)v_1 +\varphi_2 (M)v_2+\cdots+\varphi_p(M)v_p,
\end{equation}
where the matrix $M$ here could be $hJ_n$ or $c_i hJ_n$. Starting from a seminal contribution by Hochbruck and Lubich \cite{HL97} (which they analyzed Krylov subspace methods for efficiently computing the action of a matrix exponential (with a large norm) on some vector), many more efficient techniques have been proposed. A large portion of these developments is concerned with computing the expression  \eqref{eq3.14}. For example, we mention some of 
the state-of-the-art algorithms: $\mathtt{expmv}$  proposed by Al-Mohy and Higham in \cite{AH11} (using a truncated standard Taylor series expansion), $\mathtt{phipm}$ proposed by Niessen and Wright in \cite{NW12} (using adaptive Krylov subspace methods), and $\mathtt{expleja}$ proposed by Caliari et al.~in \cite{Caliari2009,CKOS16} (using Leja interpolation). With respect to computational time, it turns out that  $\mathtt{phipm}$ offer an advantage. This algorithm utilizes an adaptive time-stepping method to evaluate  \eqref{eq3.14} using only one matrix function (see Subsection~\ref{subsec:time-stepping} below). This task is carried out in a lower dimensional Krylov subspace using standard Krylov subspace projection methods i.e. the Arnoldi iteration. Moreover, the dimension of  Krylov subspaces and the number of substeps are also chosen adaptivity for improving efficiency.

Recently, the $\mathtt{phipm}$ routine was modified by Gaudreault and Pudykiewicz in \cite{GaudreaultPudykiewicz16} (Algorithm 2) by using the incomplete orthogonalization method (IOM) within the Arnoldi iteration and by adjusting the two crucial initial parameters for starting the Krylov adaptivity. This results in the new routine called $\mathtt{phipm/IOM2}$.  It is shown in \cite{GaudreaultPudykiewicz16} that this algorithm reduces computational time significantly compared to $\mathtt{phipm}$ when integrating the shallow water equations on the sphere. 

Very recently, the authors of \cite{Luan18} further improved $\mathtt{phipm/IOM2}$ which resulted in a more efficient routine named as $\mathtt{phipm\_simul\_iom2}$. For the reader's convenience, we present the idea of the adaptive time-stepping method (originally proposed in \cite{NW12}) for evaluating  \eqref{eq3.14} and introduce some new features of the new routine $\mathtt{phipm\_simul\_iom2}$.

\subsubsection{Time-stepping-based Computing of Linear $\varphi$-combinations}
\label{subsec:time-stepping}
It was observed  that the following linear ODE
\begin{equation} \label{eq3.15}
u'(t)=Mu(t)+v_1+t v_2+\cdots+\frac{t^{p-1}}{(p-1)!}v_p, \ u(0)=v_0,
\end{equation}
defined on the interval $[0, 1]$ has the exact solution at $t=1$, $u(1)$ to be the expression \eqref{eq3.14}. The time-stepping technique approximates $u(1)$ by discretizing $[0,1]$ into subintervals $0=t_0<t_1<\cdots<t_k<t_{k+1}=t_k+\tau_k<\cdots<t_K=1$ with a substepsize sequence $\tau_k$ ($k=0,1,\ldots,K-1$) and using the following relation between $u(t_{k+1})$ and its previous solution $u(t_k)$:
\begin{equation} \label{eq3.16}
u(t_{k+1})=\varphi_0 (\tau_k M)u(t_k)+\sum_{i=1}^{p}\tau^i_k \varphi_i (\tau_k M) \sum_{j=0}^{p-i} \frac{t^j_k}{j!}v_{i+j}.
\end{equation}
Using the recursion relation \eqref{eq3.7},  \eqref{eq3.16} can be simplified  as
\begin{equation} \label{eq3.17}
u(t_{k+1})=\tau^p_k \varphi_p (\tau_k M)w_p+ \sum_{j=0}^{p-i} \frac{t^j_k}{j!}w_j,
\end{equation}
where the vectors $w_j$ satisfy the recurrence relation
\begin{equation} \label{eq3.18}
w_0=u(t_k), \ w_j=Mw_{j-1}+\sum_{\ell=0}^{p-j} \frac{t^{\ell}_k}{\ell!}v_{j+\ell}, \ j=1,\ldots,p.
\end{equation}
Equation \eqref{eq3.17} implies that evaluating $u(t_K)=u(1)$ i.e. the expression \eqref{eq3.14} requires only one matrix function $\varphi_p (\tau_k A)w_p$ in each substep instead of $(p+1)$ matrix-vector multiplications. As  $0<\tau_k<1$, this task can be carried out in a Krylov subspace of lower dimension $m_k$, and in each substep only one Krylov projection is needed. With a reasonable number of substeps  $K$, it is thus expected that the total computational cost of  $\mathcal{O}(m^2_1)+\cdots+\mathcal{O}(m^2_K)$ for approximating $\varphi_p (\tau_k M)w_p$ is less than that of  $\mathcal{O}(m^2)$  for approximating $\varphi_p (M)v$  in a Krylov subspace of dimension $m$. If $K$ is too large (e.g.~ when the spectrum of $M$ is very large), this might be not true. This case, however, is handed by using the adaptive Krylov algorithm in \cite{NW12}  allowing to adjust both the dimension $m$ and the step sizes $\tau_k$ adaptivity.
This explains the computational advance of this approach compared to standard Krylov algorithms.
 
\subsubsection{New Routine $\mathtt{phipm\_simul\_iom2}$ \cite{Luan18}}
\label{sec:implementation}
 Motivated by the  two observations mentioned in items (i) and (ii) below, the $\mathtt{phipm/IOM2}$  routine (see \cite{GaudreaultPudykiewicz16}) was  modified in \cite{Luan18} for more efficient implementation of exponential Rosenbrock methods. The resulting routine \texttt{phipm\_simul\_iom2}  optimizes computational aspects of $\mathtt{phipm/IOM2}$  corresponding to these observations. In particular, along with the motivation we also recall the  two specific changes from the \texttt{phipm/IOM2}:
 \begin{itemize}
\item [(i)] Unlike \eqref{eq3.14}, where  each of the $\varphi_k$ functions is evaluated at the same argument $M$, the internal stages of exponential Rosenbrock schemes require evaluating the $\varphi$ functions at fractions of the matrix $M$:
\begin{equation}
\label{eq:exp_multistage}
   w_k = \sum_{l=1}^p \varphi_l(c_k\, M) v_l, \quad k=2,\ldots,s,
\end{equation}
where now the node values $c_2, \ldots, c_s$  are scaling factors used for each $v_k$ output. 
To optimize this evaluation, \texttt{phipm\_simul\_iom2} computes all $w_k$ outputs in  \eqref{eq:exp_multistage} simultaneously, instead of computing only one at a time. This is accomplished by first requiring that the entire array $c_2, \ldots, c_s$ as an input to the function.  
Within the substepping process \eqref{eq3.16}, each value $c_j$  is  aligned with a substep-size $\tau_k$.  The solution vector is stored at each of these moments and on output the full set $\{w_k\}_{k=1}^s$ is returned.  Note that this approach is similar but differs from \cite{TLP12} that it guarantees no loss of solution accuracy since it explicitly stops at each $c_k$ instead of using interpolation to compute  $w_k$ as in \cite{TLP12}.
\item [(ii)]~In view of the higher-order exponential Rosenbrock schemes (see also from Section
\ref{sec:selected_exponential_schemes}), it is realized that they usually use a subset of the
$\varphi_l$ functions. Therefore,  multiple vectors in \eqref{eq:exp_multistage} will be zero.
In this case, \texttt{phipm\_simul\_iom2} will check whether $w_{j-1}\ne 0$  (within the recursion \eqref{eq3.18}) before computing the matrix-vector product $M\,w_{j-1}$.  While matrix-vector products require $\mathcal O(n^2)$ work,  checking $u\ne 0$ requires only $\mathcal O(n)$. This can result in significant savings for large $n$.
 \end{itemize}

\subsubsection{Implementation of $\mathtt{exprb42}$ and $\mathtt{pexprb43}$} 
Taking a closer look at the structures of the two selected exponential Rosenbrock schemes $\mathtt{exprb42}$ and $\mathtt{pexprb43}$, we now
 make use of \texttt{phipm\_simul\_iom2} for implementing these schemes. For simplicity, let us denote  $M= h J_n$ and $v=hF(u_n)$.\\
\\
\emph{Implementation of $\mathtt{exprb42}$:}
Due to the structure of $\mathtt{exprb42}$ given in  \eqref{eq3.11}, one needs two calls to \texttt{phipm\_simul\_iom2}:
 \begin{itemize}
\item [(i)] Evaluate  $y_1=\varphi _{1} (\tfrac{3}{4} M)w_1$ with $w_1= \tfrac{3}{4} v$ (so $w_0=0$) to get  $U_{n2}=u_n+ y_1$,
\item [(ii)]~Evaluate  $w=\varphi _{1} (M)v_1+ \varphi _{3} (M)v_3$  (i.e. $v_0=v_2=0$) with  $v_1=v, v_3=\tfrac{32}{9} hD_{n2}$ to get  $u_{n+1}=u_n+ w.$
 \end{itemize}
\emph{Implementation of $\mathtt{pexprb43}$:}
Although $\mathtt{pexprb43}$ is a 3-stage scheme, its special structure  \eqref{eq3.12} allows to use only two calls to \texttt{phipm\_simul\_iom2}:
 \begin{itemize}
\item [(i)]~Evaluate both terms $y_1=\varphi _{1} (\tfrac{1}{2} M)v$ and $z_1=\varphi _{1} (M)v$ simultaneously  to get the two stages  $U_{n2}=u_n+\tfrac{1}{2}y_1$ and $U_{n3}=u_n+z_1$,
\item [(ii)]~Evaluate  $w=\varphi _{3} (M)v_3+ \varphi _{4} (M)v_4$   (i.e. $v_0=v_1=v_2=0$) with $v_3=h(16D_{n2}-2D_{n3}), v_4=h(-48D_{n2}+12D_{n3})$ to get $u_{n+1}=U_{n3}+w. $
 \end{itemize}

\section{Numerical Examples}
\label{sec:experiments}

In this section we present a broad spectrum of numerical examples to study the behavior of the presented exponential Rosenbrock-type methods, in particular the fourth-order scheme $\mathtt{exprb42}$ using two stages and the fourth-order $\mathtt{pexprb43}$ scheme using three stages implemented in parallel.

In particular, we focus on relevant aspects in the realm of visual computing, like stability and energy conservation, large stiffness, and high fidelity and visual accuracy. A tabular summary of the models that are used throughout this section can be found in Table~\ref{tab:TestCases}. Furthermore, our simulation includes important aspects like elastic collisions and nonelastic deformations. The presented exponential Rosenbrock-type methods are evaluated against classical and state-of-the-art methods used in visual computing, in particular against the implicit-explicit variational (IMEX) integrator (cf. \cite{Stern:2006,Stern:2009}), the standard fourth-order Runge--Kutta method (see \cite{Runge:1895,Kutta:1901}), and the implicit BDF-1 integrator (see \cite{Curtiss:1952}). All simulation results visualized here have been computed using a machine with an Intel(R) Xeon E5 3.5 GHz and 32 GB DDR-RAM. For each simulation scenario the largest possible time step size is used which still leads to a desired visually plausible result.

\begin{table}[H]
\caption{Overview of the test cases used for the numerical experiments. Their complexity $N$ (i.e.~the number of the resulting equations of motion), the simulated time, and respective running times for the exponential Rosenbrock-type methods $\mathtt{exprb42}$ and $\mathtt{pexprb43}$, the implicit-explicit variational integrator (V-IMEX), the standard fourth order Runge--Kutta method (RK4), and the BDF-1 integrator are shown.}
\begin{center}
\begin{tabular}{ |c|c|c|c|c|c|c|c|c| }
\hline
\textbf{No.} & \textbf{Model} & $N$ & \textbf{Sim. Time} & $\mathtt{exprb42}$  & $\mathtt{pexprb43}$ & \textbf{V-IMEX} & \textbf{RK 4} & \textbf{BDF-1} \\ 
\hline
1 & Coil Spring 	& 24k	& 60\,s	& 55\,s	& 47\,s	& 12\,min	& 46\,min	& 62\,min \\
\hline
2 & Brushing	& 90k	& 15\,s	& 52\,s	& 51\,s	& 11\,min	& 53\,min	& 72\,min \\
\hline
3 & Crash Test (moderate)	& 360k	& 2\,s	& 44\,s	& 44\,s	& 9\,min	& 47\,min 	& 58\,min \\
4 & Crash Test (fast)	& 360k	& 2\,s	& 47\,s	& 46\,s	& 9\,min	& 47\,min	& 59\,min \\
\hline
\end{tabular}
\end{center}
\label{tab:TestCases}
\end{table}%

\subsection{Simulation of Deformable Bodies}

\begin{figure*}
\centering
\includegraphics[width=0.95\textwidth]{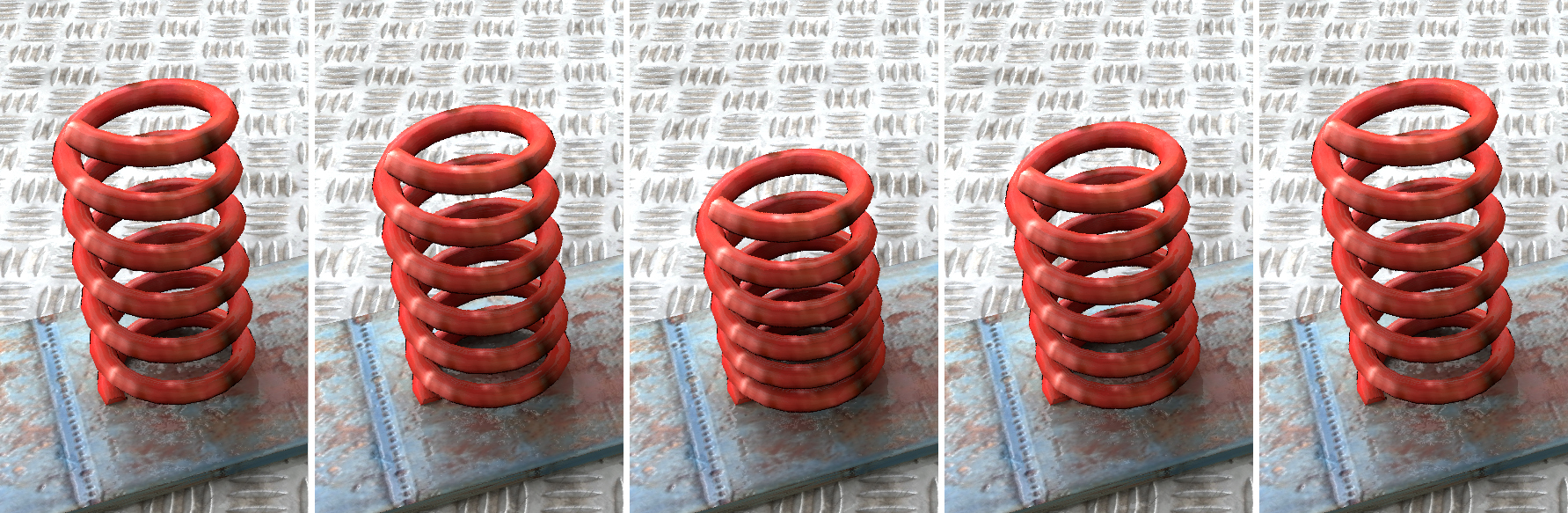}
\caption{Simulation of an oscillating coil spring.}
\label{fig:CoilSpring}
\end{figure*}

In order to illustrate the accurate energy preservation of the presented exponential Rosenbrock-type methods, we set up an undamped scene of an oscillating coil spring, which is modeled as a deformable body composed of tetrahedra, in particular of $8\,000$ vertices corresponding to $N=24\,000$ equations of motion, which are derived from a system of coupled oscillators with uniform spring stiffness of $k=10^6$. Since the coil spring is exposed to an external forces field, it starts to oscillate as illustrated in Figure~\ref{fig:CoilSpring}. It can be seen that the top of the coil spring returns to its initial height periodically during the simulation which can be seen as an indicator for energy conservation. In fact when using the exponential Rosenbrock-type methods $\mathtt{exprb42}$ and $\mathtt{pexprb43}$ we observe that the discrete energy is only slightly oscillating around the real energy without increasing oscillations over time. In contrast, the standard fourth order Runge--Kutta method respectively the BDF-1 integrator generate significant numerical viscosity leading to a loss of energy around $22\%$ respectively $40\%$ after 60\,s of simulated time.

The exponential Rosenbrock-type methods $\mathtt{exprb42}$ and $\mathtt{pexprb43}$ show their advantageous behavior since these methods can be applied with orders of magnitude larger time steps compared to the other integrators. Even with a step size of $h=0.05$ the relative error is still below $2\%$ for $\mathtt{exprb42}$ and about a single percent for $\mathtt{pexprb43}$.\footnote{We estimated the error after 60\,s of simulated time based on the accumulated Euclidean distances of the individual particles in the position space compared to ground truth values which are computed with a sufficiently small step size.} From a point of view of computation time, we achieve a speed up of a factor of around thirteen using $\mathtt{exprb42}$ and of over fifteen using $\mathtt{pexprb43}$ compared to the second best method, the variational IMEX integrator as illustrated in Table~\ref{tab:TestCases}. Compared to the other methods, the exponential Rosenbrock-type methods allow for accurate simulations in real-time.

\subsection{Simulation of Fibers including Elastic Collisions}

Fibers are canonical examples for complex interacting systems. According to the work of Michels et al.~(see \cite{Michels:2015b}), we set up a toothbrush composed of individual bristles. Each bristle consists of coupled oscillators that are connected in such a way that the fiber axis is enveloped by a chain of cuboid elements. For preventing a volumetric collapse during the simulation, additional diagonal springs are used. The toothbrush consists of $1\,500$ bristles, each of $20$ particles leading to $90\,000$ equations of motion. We make use of additional repulsive springs in order to prevent from interpenetrations.\footnote{In order to detect collisions efficiently, we make use of a standard bounding volume hierarchy.} Since the approach allows for the direct use of realistic parameters in order to set up the stiffness values in the system of coupled oscillators, we employ a Young's modulus of $3.2\cdot10^6\,\text{Ncm$^{-2}$}$, a torsional modulus of $10^5\,\text{Ncm$^{-2}$}$, and segment thicknesses of 0.12\,mm. 

\begin{figure*}
\centering
\includegraphics[width=0.95\textwidth]{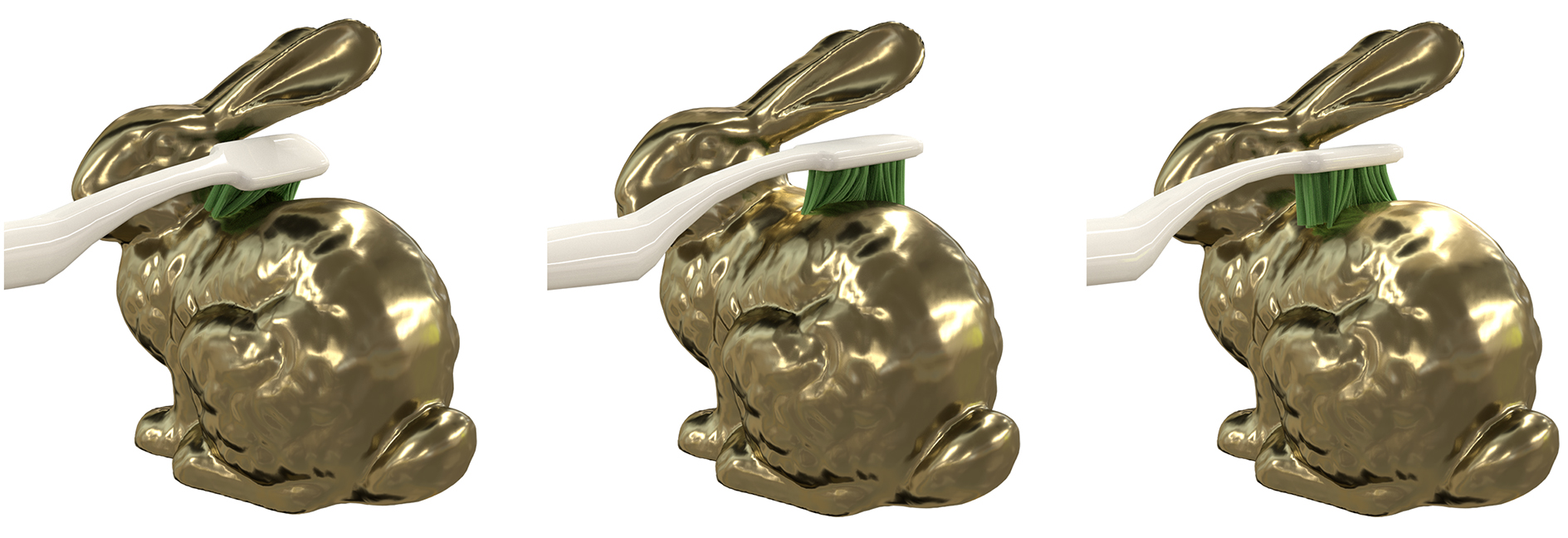}
\caption{Simulation of a brush cleaning a bronze-colored paperweight.}
\label{fig:BunnyBrushing}
\end{figure*}

We simulate 15\,s of a toothbrush cleaning a paperweight illustrated in Figure~\ref{fig:BunnyBrushing}. This simulation can be carried out almost in real-time which is not possible with the use of classical methods as illustrated in Table~\ref{tab:TestCases}.

\subsection{Crash Test Simulation including Nonelastic Deformations}

\begin{figure*}
\centering
\includegraphics[width=0.95\textwidth]{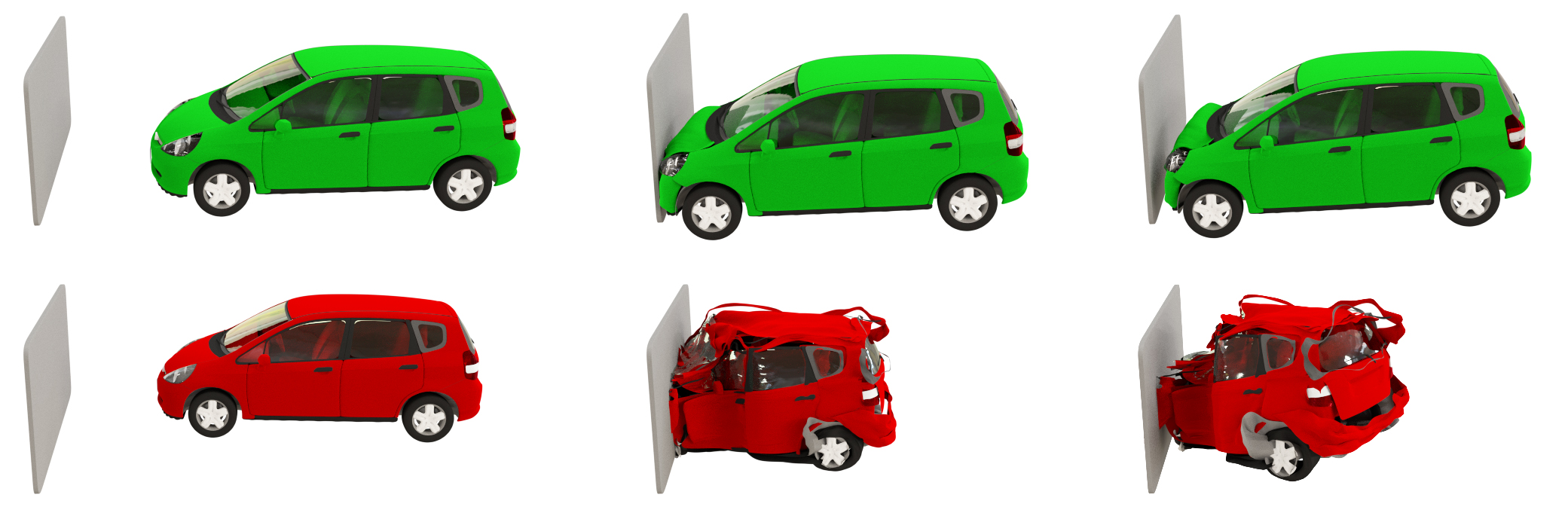}
\caption{Simulations of two frontal nonelastic crash scenarios: a car with moderate velocity (top) and high velocity (bottom).}
\label{fig:CrashTests}
\end{figure*}

As a very complex example with relevance in the context of special effects, we simulate a frontal crash of a car into a wall as illustrated in Figure~\ref{fig:CrashTests}. The mesh of the car and its interior is composed of $120\,000$ vertices leading to $360\,000$ equations of motion. The global motion (i.e.~the rebound of the car) is computed by treating the car as a rigid body. Using an appropriate bounding box, this can be easily carried out in real-time. The deformation is then computed using a system of coupled oscillators with structural stiffness values of $k=10^4$ and bending stiffness values of $k/100$. If the deformation reaches a defined threshold, the rest lengths of the corresponding springs are corrected in a way, that they do not elastically return to their initial shape. Using the exponential Rosenbrock-type methods, the whole simulation can be carried out at interactive frame rates. Such an efficient computation can not be achieved with established methods as illustrated in Table~\ref{tab:TestCases}.

\section{Conclusion}
\label{sec:conclusion}

We introduced  the class of \textit{explicit exponential Rosenbrock methods} for the time integration of large systems of nonlinear differential equations. In particular, the exponential Rosenbrock-type fourth-order schemes $\mathtt{exprb42}$ using two stages and $\mathtt{pexprb43}$ using three stages were discussed and their implementation were addressed. In order to study their behavior, a broad spectrum of numerical examples was computed. In this regard, the simulation of deformable bodies, fibers including elastic collisions, and crash scenarios including nonelastic deformations was addressed focusing on relevant aspects in the realm of visual computing, like stability and energy conservation, large stiffness values, and high fidelity and visual accuracy. An evaluation against classical and state-of-the-art methods was presented demonstrating their superior performance with respect to the simulation of large systems of stiff differential equations.

\bibliographystyle{spmpsci}
\bibliography{references}
\end{document}